\documentclass[12pt,leqno]{article}
\usepackage{amssymb}
\usepackage[english,francais]{babel}
\newtheorem{theorem}{Theorem}[section]
\newtheorem{lemma}[theorem]{Lemma}
\newtheorem{e-proposition}[theorem]{Proposition}

\newtheorem{e-definition}[theorem]{Definition\rm}

\newtheorem{theoreme}{Th\'eor\`eme}[section]

\newtheorem{proposition}[theoreme]{Proposition}

\def\og{\leavevmode\raise.3ex\hbox{$\scriptscriptstyle\langle\!\langle$~}}
\def\fg{\leavevmode\raise.3ex\hbox{~$\!\scriptscriptstyle\,\rangle\!\rangle$}}
\begin{document}
\title{Error calculus and regularity of Poisson functionals:\\
the lent particle method.}
\author{Nicolas BOULEAU\footnote{Ecole des Ponts, Paris-Est, ParisTech. email: {\tt bouleau@enpc.fr}}}
\date{--}
\selectlanguage{english}
\maketitle

{\vskip 0.5\baselineskip
\noindent{\bf Abstract} \vskip 0.5\baselineskip \noindent
\selectlanguage{english}
We propose a new method to apply the Lipschitz functional calculus of local Dirichlet forms to Poisson random measures.

\vskip 0.5\baselineskip

\selectlanguage{francais}
\noindent{\bf R\'esum\'e} \vskip 0.5\baselineskip \noindent
{\bf Calcul d'erreur et r\'egularit\'e des fonctionnelles de Poisson: la m\'ethode de la particule pr\^et\'ee. }
Nous proposons une nouvelle m\'ethode pour appliquer le calcul fonctionnel lipschitzien des formes de Dirichlet locales aux mesures al\'eatoires de Poisson.
}


\selectlanguage{english}
\section{Notation and basic formulae.}
\label{}
Let us consider a local Dirichlet structure with carr\'e du champ $(X,\mathcal{X},\nu,\mathbf{d},\gamma)$ where $(X,\mathcal{X},\nu)$ is a $\sigma$-finite measured space called {\it bottom-space}. Singletons are in $\mathcal{X}$ and $\nu$ is diffuse, $\mathbf{ d}$ is the domain of the Dirichlet form $\epsilon[u]=1/2\int \gamma[u]d\nu$. We denote $(a,\mathcal{D}(a))$ the generator in $L^2(\nu)$ (cf. \cite{B-H}).

A random Poisson measure associated to $(X,\mathcal{X},\nu)$ is denoted $N$. $\Omega$ is the configuration space of countable sums of Dirac masses on $X$ and $\mathcal{A}$ is the $\sigma$-field generated by $N$, of law $\mathbb{P}$ on $\Omega$. The space $(\Omega, \mathcal{A},\mathbb{P})$ is called {\it the up-space}. We write $N(f)$ for $\int f dN$.
If $p\in[1,\infty[$ the set $\{e^{iN(f)} : f {\mbox{ real, }} f\in L^1\cap L^2(\nu)\}$ is total in $L^p_\mathbb{C}(\Omega,\mathcal{A},\mathbb{P})$. We put $\tilde{N}=N-\nu$. The relation $\mathbb{E}(\tilde{N}f)^2=\int f^2d\nu$ extends and gives sense to $\tilde{N}(f)$, $f\in L^2(\nu)$. The Laplace functional and the differential calculus with $\gamma$ yield

\begin{equation} \forall f\in\mathbf{ d}, \forall h\in \mathcal{D}(a)
\quad\quad
\mathbb{E}[e^{i\tilde{N}(f)}(\tilde{N}(a[h])+\frac{i}{2}N(\gamma[f,h])]=0.
\end{equation}

\section{Product, particle by particle, of a Poisson random measure by a probability measure.}
 Given a probability space  $(R,\mathcal{R},\rho)$, let us consider a Poisson random measure $N\odot\rho$ on $(X\times R,\mathcal{X}\times\mathcal{R})$ with intensity $\nu\times\rho$ such that for $f\in L^1(\nu)$ and $g\in L^1(\rho)$ if $N(f)=\sum f(x_n)$ then $(N\odot\rho)(fg)=\sum f(x_n)g(r_n)$ where the $r_n$'s are i.i.d.  independent of $N$ with law $\rho$. Calling $(\hat{\Omega},\hat{\mathcal{A}}, \hat{\mathbb{P}})$ the product of all the factors $(R,\mathcal{R},\rho)$ involved in the construction of $N\odot\rho$, we obtain the following properties :
 For an $\mathcal{A}\times\mathcal{X}\times\mathcal{R}$-measurable and positive function $F$, 
$
\hat{\mathbb{E}}\int F(\omega, x, r) N\odot\rho(dxdr)=\int F\;d\rho\,dN\quad\mathbb{P}\mbox{-a.s.}
$

Let us denote by $\mathbb{P}_N$ the measure $\mathbb{P}(d\omega)N_\omega(dx)$ on $(\Omega\times X, \mathcal{A}\times\mathcal{X})$. We have
\begin{lemma} Let $F$ be $\mathcal{A}\times\mathcal{X}\times\mathcal{R}$-measurable, $F\in\mathcal{L}^2(\mathbb{P}_N\times\rho)$ and such that 

$\int F(\omega, x, r)\;\rho(dr)=0\quad \mathbb{P}_N$-a.s., then $\int F\;d(N\odot\rho)$ is well defined, belongs to $L^2(\mathbb{P}\times\hat{\mathbb{P}})$ and
\begin{equation}
\hat{\mathbb{E}}(\int F\;d(N\odot\rho))^2=\int F^2\;dN\,d\rho\qquad\mathbb{P}\mbox{-a.s.}
\end{equation}
\end{lemma}
The argument consists in considering $F_n$ satisfying 

$\mathbb{E}\int F_n^2\;d\nu d\rho<+\infty$ and $\mathbb{E}\int (\int |F_n|\;d\nu)^2 d\rho<+\infty$ and then using the relation

$\hat{\mathbb{E}}(\int F_n\;d(N\odot\rho))^2=(\int F_n d\rho dN)^2-\int(\int F_n d\rho)^2dN+\int F_n^2 d\rho dN\quad \mathbb{P}\mbox{-a.s.}$

\section{Construction by Friedrichs' method and expression of the gradient.} 
a) We suppose the space by $\mathbf{ d}$ of the bottom structure is separable, then a gradient exists (cf. \cite{B-H} Chap. V, p.225 {\it et seq}.). We denote it $\flat$ and choose it with values in the space $L^2(R,\mathcal{R},\rho)$. Thus, for $u\in \mathbf{ d}$ we have $u^\flat\in L^2(\nu\times\rho)$, $\gamma[u]=\int(u^\flat)^2d\rho$  and $\flat$ satisfies the chain rule. We suppose in addition, what is always possible, that $\flat$ takes its values in the subspace orthogonal to the constant $1$, i.e.
\begin{equation}
\forall u\in \mathbf{ d}\qquad \int u^\flat\;d\rho=0\quad \nu\mbox{-a.s.}
\end{equation}
This hypothesis is important here as  in many applications (cf. \cite{B} Chap V \S 4.6). We suppose also, but this is not essential (cf. \cite{B-H} p44)
$
1\in \mathbf{ d}_{loc}\quad \gamma[1]=0 \mbox{ so that } 1^\flat=0
.$

b) We define a pre-domain $\mathcal{D}_0$ dense in $L^2_\mathbb{C}(\mathbb{P})$ by
$$\mathcal{D}_0=\{\sum_{p=1}^m \lambda_pe^{i\tilde{N}(f_p)} ; m\in\mathbb{N}^\ast, \lambda_p\in\mathbb{C}, f_p\in \mathcal{D}(a)\cap L^1(\nu)\}.$$

c) We introduce the creation operator inspired from quantum mechanics (see \cite{Nualart}, \cite{Picard}, \cite{Privault}, \cite{Albeverio}, \cite{Ishikawa},\cite{Ma} and \cite{Sole} among others) defined as follows
\begin{equation}
\varepsilon^+_x(\omega)\mbox{ equals }\omega \mbox{ if }x\in\mbox{supp}(\omega),\;\mbox{ and equals }\omega+\varepsilon_x\mbox{ if }x\notin\mbox{supp}(\omega)
\end{equation}
so that 
\begin{equation}
\quad \varepsilon^+_x(\omega)=\omega\quad N_\omega\mbox{-a.e. } x
\quad\mbox{ and }\quad
\varepsilon^+_x(\omega)=\omega+\varepsilon_x\quad\nu\mbox{-a.e. } x
\end{equation}
This map is measurable and the Laplace functional shows that for an $\mathcal{A}\times\mathcal{X}$-measurable $H\geq 0$, 
\begin{equation}
\mathbb{E}\int \varepsilon^+H\;d\nu=\mathbb{E}\int H\;dN.
\end{equation} Let us remark also that by (5), for $F\in \mathcal{L}^2(\mathbb{P}_N\times\rho)$
\begin{equation}
\int \varepsilon^+ F\;d(N\odot\rho)=\int Fd(N\odot\rho)\qquad \mathbb{P}\times\hat{\mathbb{P}}\mbox{-a.s.}
\end{equation}

d) We defined a gradient $\sharp$ for the up-structure on $\mathcal{D}_0$ by
putting for $F\in\mathcal{D}_0$
\begin{equation}
F^\sharp=\int (\varepsilon^+F)^\flat\;d(N\odot\rho)
\end{equation}
this definition being justified by the fact that for $\mathbb{P}$-a.e. $\omega$ the map $y\mapsto F(\varepsilon^+_y(\omega))-F(\omega)$ is in $\mathbf{ d}$, $\varepsilon^+F$ belongs to $L^\infty(\mathbb{P})\otimes \mathbf{ d}$ algebraic tensor product, and $(\varepsilon^+F-F)^\flat=(\varepsilon^+F)^\flat\in L^2(\mathbb{P}_N\times\rho)$.

For $F,G\in\mathcal{D}_0$ of the form
$$F=\sum_p \lambda_pe^{i\tilde{N}(f_p)}=\Phi(\tilde{N}(f_1),\ldots, \tilde{N}(f_m))\qquad G=\sum_q \mu_qe^{i\tilde{N}(g_q)}=\Psi(\tilde{N}(g_1),\ldots, \tilde{N}(g_n))$$ we compute using (2), (3) and (7) (in the spirit of prop. 1 of \cite{Privault} or lemma 1.2 of \cite{Ma})
\begin{equation}
\hat{\mathbb{E}}[F^\sharp\overline{G^\sharp}]=\sum_{p,q}\lambda_p\overline{\mu_q}e^{i\tilde{N}(f_p)-i\tilde{N}(g_q)}\gamma[f_p,g_q]
\end{equation} and we have
\begin{proposition} If we put
$A_0[F]=\sum_p\lambda_pe^{i\tilde{N}(f_p)}(i\tilde{N}(a[f_p])-\frac{1}{2}N(\gamma[f_p]))$ it comes
\begin{equation}
\mathbb{E}[A_0[F]\overline{G}]=-\frac{1}{2}\mathbb{E}\sum_{p,q}\Phi^\prime_p\overline{\Psi^\prime_q}N(\gamma[f_p,g_q]).
\end{equation}
\end{proposition}
In order to show that $A_0[F]$ does not depend on the form of $F$, by (10) it is enough to show that the expression $\sum_{p,q}\Phi^\prime_p\overline{\Psi^\prime_q}N(\gamma[f_p,g_q])$  depends only on $F$ and $G$. But this comes from (9) since $F^\sharp$ and $G^\sharp$ depend only on $F$ and $G$.

By this proposition, $A_0$ is symmetric on $\mathcal{D}_0$, negative, and the argument of Friedrichs applies (cf \cite{B-H} p4), $A_0$ extends uniquely to a selfadjoint operator $(A,\mathcal{D}(A))$ which defines a closed positive (hermitian) quadratic form $\mathcal{E}[F]=-\mathbb{E}[A[F]\overline{F}]$. By (10) contractions operate and (cf. \cite{B-H}) $\mathcal{E}$ is a Dirichlet form which is local with carr\'e du champ denoted $\Gamma$ and the up-structure obtained $(\Omega, \mathcal{A}, \mathbb{P}, \mathbb{D}, \Gamma)$ satisfies
\begin{equation}
\forall f\in \mathbf{ d}, \quad \tilde{N}(f)\in\mathbb{D}\;\mbox{ and } \Gamma[\tilde{N}(f)]=N(\gamma[f])
\end{equation} The operator $\sharp$ extends to a gradient for $\Gamma$ as a closed operator from $L^2(\mathbb{P})$ into $L^2(\mathbb{P}\times\hat{\mathbb{P}})$ with domain $\mathbb{D}$ which satisfies the chain rule and may be computed on functionals $\Phi(\tilde{N}(f_1),\ldots,\tilde{N}(f_m))$, $\Phi$ Lipschitz and $\mathcal{C}^1$ and their limits in $\mathbb{D}$ (as done in \cite{Denis}).

Formula (8) for $\sharp$ can be extended from $\mathcal{D}_0$ to $\mathbb{D}$. Let us introduce the space $\underline{\mathbb{D}}$  closure of $\mathcal{D}_0\otimes\mathbf{ d}$ for the norm
$$\|H\|_{\underline{\mathbb{D}}}
=(\mathbb{E}\int\gamma[H(\omega,.)](x)\;N(dx))^{1/2}+\mathbb{E}\int |H(\omega,x)|\xi(x)\;N(dx)$$ where $\xi>0$ is a fixed function such that $N(\xi)\in L^2(\mathbb{P})$.\\
\begin{theorem} The formula 
$F^\sharp=\int(\varepsilon^+F)^\flat\;d(N\odot\rho)$ decomposes as follows
$$F\in\mathbb{D}\quad\stackrel{\varepsilon^+}{\longmapsto}\quad \varepsilon^+F\in\underline{\mathbb{D}}\quad\stackrel{\flat}{\longmapsto}\quad(\varepsilon^+F)^\flat\in L^2_0(\mathbb{P}_N\times\rho)\quad\stackrel{d(N\odot\rho)}{\longmapsto}\quad F^\sharp\in L^2(\mathbb{P}\times\hat{\mathbb{P}})$$
where each operator is continuous on the range of the preceding one, $L^2_0(\mathbb{P}_N\times\rho)$ denoting the closed subspace of $L^2(\mathbb{P}_N\times\rho)$  of $\rho$-centered elements, and we have
\begin{equation}
\Gamma[F]=\hat{\mathbb{E}}|F^\sharp|^2=\int\gamma[\varepsilon^+F]\;dN.
\end{equation}
\end{theorem}
\section{The lent particle method.}
Let us consider, for instance, a real process $Y_t$ with independent increments and L\'evy measure $\sigma$ integrating $x^2$, $Y_t$ being supposed centered without Gaussian part. We assume that $\sigma$ has an l.s.c. density so that a local Dirichlet structure may be constructed on $\mathbb{R}\backslash\{0\}$ with carr\'e du champ 
$
\gamma[f]=x^2f^{\prime 2}(x).
$
If $N$ is the random Poisson measure with intensity $dt\times\sigma$ we have $\int_0^th(s)\;dY_s=\int1_{[0,t]}(s)h(s)x\tilde{N}(dsdx)$ and the choice done for $\gamma$ gives 
$\Gamma[\int_0^th(s)dY_s]=\int_0^th^2(s)d[Y,Y]_s
$ for $h\in L^2_{loc}(dt)$.
In order to study the regularity of the random variable $ V=\int_0^t\varphi(Y_{s-})dY_s$ where $\varphi$ is Lipschitz and $\mathcal{C}^1$, we have two ways:

a) We may  represent the gradient $\sharp$ as $Y_t^\sharp=B_{[Y,Y]_t}$ where $B$ is a standard auxiliary independent Brownian motion. Then by the chain rule
$
V^\sharp=\int_0^t\varphi^\prime(Y_{s-})(Y_{s-})^\sharp dY_s+\int_0^t\varphi(Y_{s-})dB_{[Y]_s}$ now, using $(Y_{s-})^\sharp=(Y_s^\sharp)_-$,  a classical but rather tedious stochastic computation yields
\begin{equation}
\Gamma[V]=\textstyle{\hat\mathbb{E}[V^{\sharp 2}]=\sum_{\alpha\leq t}\Delta Y_\alpha^2(\int_{]\alpha}^t\varphi^\prime(Y_{s-})dY_s+\varphi(Y_{\alpha-}))^2}.
\end{equation} Since $V$ has real values the {\it energy image density property} holds, and $V$ has a density as soon as $\Gamma[V]$ is strictly positive a.s. what may be discussed using the relation (13).

b) Another more direct way consists in applying the theorem. For this we define $\flat$ by choosing $\eta$ such that $\int_0^1\eta(r)dr=0$ and $\int_0^1\eta^2(r)dr=1$ and putting
$f^\flat=xf^\prime(x)\eta(r).
$

$1^o$. First step. We add a particle $(\alpha,x)$ i.e. a jump to $Y$ at time $\alpha$ with size $x$ what gives

\noindent$\varepsilon^+V-V=\varphi(Y_{\alpha-})x+\int_{]\alpha}^t(\varphi(Y_{s-}+x)-\varphi(Y_{s-}))dY_s
$

$2^o$. $V^\flat=0$ since $V$ does not depend on $x$, and

\noindent$(\varepsilon^+V)^\flat=(\varphi(Y_{\alpha-})x+\int_{]\alpha}^t\varphi^\prime(Y_{s-}+x)xdY_s)\eta(r)\quad$ 
\noindent because $x^\flat=x\eta(r)$.

$3^o$. We compute 
$\gamma[\varepsilon^+V]=\int(\varepsilon^+V)^{\flat2}dr=(\varphi(Y_{\alpha-})x+\int_{]\alpha}^t\varphi^\prime(Y_{s-}+x)xdY_s)^2$

$4^o$. We take back the particle we gave, because in order to compute $\int\gamma[\varepsilon^+V]dN$ the integral in $N$ confuses $\varepsilon^+\omega$ and $\omega$. 

That gives
$\int\gamma[\varepsilon^+V]dN=\int(\varphi(Y_{\alpha-})+\int_{]\alpha}^t\varphi^\prime(Y_{s-})dY_s)^2x^2\;N(d\alpha dx)$ and (13).

We remark that both operators $F\mapsto \varepsilon^+F$, $F\mapsto(\varepsilon^+F)^\flat$ are non-local, but instead $F\mapsto \int (\varepsilon^+F)^\flat\,d(N\odot\rho)$ and $F\mapsto \int\gamma[\varepsilon^+F]\,dN$ are local : taking back the lent particle gives the locality.

\end{document}